\documentclass[11pt,leqno]{amsart}
\usepackage{amsmath,amsfonts,amsthm,amssymb}
\usepackage{xspace}
\usepackage{paralist}
\usepackage{setspace}
\usepackage{times}
\usepackage{graphicx}
\usepackage{float}
\usepackage{caption}
\usepackage{subfig}
\usepackage{hyperref}

\swapnumbers

\theoremstyle{plain}

\newtheorem{sstheorem}[subsection]{Theorem}

\theoremstyle{definition}

\newtheorem{ssdefinition}[subsection]{Definition}

\theoremstyle{remark}

\newtheorem{ssremarks}[subsection]{Remarks}

\newcommand{\cF}{\ensuremath{\mathcal{F}}}
\newcommand{\cG}{\ensuremath{\mathcal{G}}}

\newcommand{\cK}{\ensuremath{\mathcal{K}}}

\newcommand{\cP}{\ensuremath{\mathcal{P}}}
\newcommand{\cQ}{\ensuremath{\mathcal{Q}}}

\newcommand{\cT}{\ensuremath{\mathcal{T}}}

\newcommand{\bbI}{\ensuremath{\mathbb{I}}\xspace}
\newcommand{\bbM}{\ensuremath{\mathbb{M}}\xspace}

\newcommand{\bbR}{\ensuremath{\mathbb{R}}\xspace}

\newcommand{\bdry}{\partial}

\newcommand{\dist}{\operatorname{dist}}

\newcommand{\spt}{\operatorname{spt}}
\newcommand{\ra}{\rightarrow}

\newcommand\elbow
   {\kern.2em\vrule width .06 em height 1.75 ex depth .20 ex
  \vrule height -.08 ex width .375 em depth .20 ex \kern .200 em}

\renewcommand{\today}{January 30, 2017}

\begin{document}

\date{\today}

\title[Explicit Determination of \((N-1)\)-Dimensional Area Minimizing Surfaces]%
{Explicit Determination in \(\bbR^{N}\) of  \((N-1)\)-Dimensional\\
Area Minimizing Surfaces with Arbitrary Boundaries}
\author{Harold R. Parks}
\address{Department of Mathematics\\
         Oregon State University\\ 
         Corvallis, OR 97331}
\email{parksh@science.oregonstate.edu}
\author{Jon T. Pitts}
\address{Department of Mathematics\\
         Texas A\&M University\\
         College Station, TX 77843}
\email{jpitts@math.tamu.edu}


\thanks{The work of the second author was supported in part by a grant
from the National Science Foundation.}

\begin{abstract}
%

  Let \(N\ge3\)
  be an integer and \(B\)
  be a smooth, compact, oriented, \((N-2)\)-dimensional
  boundary in \(\bbR^{N}\).
  In 1960, H.~Federer and W.~Fleming \cite{FedererFleming60} proved
  that there is an \((N-1)\)-dimensional
  integral current spanning surface of least area.  The proof was by
  compactness methods and non-constructive.  In 1970 H.~Federer
  \cite{Federer70} proved the definitive regularity result for such a
  codimension one minimizing surface.  Thus it is a question of long
  standing whether there is a numerical algorithm that will closely
  approximate the area minimizing surface.  The principal result of
  this paper is an algorithm that solves this problem.
  
  Specifically, given a neighborhood~\(U\) around \(B\)
  in~\(\bbR^{N}\) and a tolerance~\(\epsilon>0\),
  we prove that one can explicitly compute in finite time an
  \((N-1)\)-dimensional integral current \(T\)
  with the following approximation requirements:
  \begin{enumerate}
  \item \(\spt(\bdry T)\subset U\).
  \item \(B\) and \(\bdry T\) are within distance \(\epsilon\) in the
     Hausdorff distance.
  \item \(B\) and \(\bdry T\) are within distance \(\epsilon\) in the
     flat norm distance.
  \item \(\bbM(T)<\epsilon+\inf\{\bbM(S):\bdry S=B\}\).
  \item Every area minimizing current \(R\)
    with \(\partial R=\partial T\)
    is within flat norm distance \(\epsilon\) of~\(T\).
  \end{enumerate}
\end{abstract}

\maketitle


\section{Introduction}\label{sec:intro}

In this paper, we will follow the notation and terminology of Federer
\cite{Federer69} except as otherwise noted.  Fix a positive
integer~\(N\geq3\).
In 1960, H.~Federer and W.~Fleming \cite{FedererFleming60} proved that
for any smooth, compact, \((N-2)\)-dimensional,
oriented boundary in \(\bbR^{N}\),
there is an \((N-1)\)-dimensional
spanning surface of least area.  The proof was by compactness methods
and non-constructive.  In 1970 H.~Federer \cite{Federer70} 
proved the definitive regularity result for such a codimension one
minimizing surface.  Thus it is a question of long standing whether
there is a numerical algorithm that will closely approximate the area
minimizing surface.  The principal result of this paper is an
algorithm that solves this problem:

\begin{sstheorem}[Main Result]
  \label{thm:main}
  Given a smooth \((N-2)\)-dimensional
  integral boundary \(B\), neighborhood \(U\)
  around \(B\), and \(\epsilon>0\),
  we will compute in finite time an integral current \(T\)
  that we can guarantee satisfies the following requirements:
\begin{enumerate}
\item\label{thm:main:1} \(\spt(\bdry T)\subset U\).
\item\label{thm:main:2} 
  \(\dist_{H}[\spt(\bdry T),\spt(B)]<\epsilon\), where
  \(\dist_{H}[\cdot,\cdot]\) is Hausdorff distance.
\item\label{thm:main:3} \(\cF(\bdry T-B)<\epsilon\).
\item\label{thm:main:4}
  \(\bbM(T)<\epsilon+\inf\{\bbM(S):\bdry S=B\}\).
\item\label{thm:main:5} Every area minimizing current \(R\)
    with \(\partial R=\partial T\)
    is within flat norm distance
    \(\epsilon\) of~\(T\).
\end{enumerate}
\end{sstheorem}

\begin{ssremarks}
  We should note that there is a limit to what can be expected.  For
  general boundary curves, the best reasonable result is the
  approximation, in both area and location, of an area minimizing
  surface that has boundary near the given boundary and has area
  nearly equal to the minimum of areas of surfaces spanning the given
  boundary.

  \begin{itemize}
  \item In general, there will be little \textit{a~priori} control of
    the topology of a minimizing surface.
  \item In general, the area minimizing surface with a given boundary
    is not unique.  Even though F.~Morgan \cite{morgan81} has shown
    that for a generic boundary the area minimizing surface is unique,
    there are but few situations in which uniqueness can be guaranteed
    \textit{a~priori}.
  \item Distinct small perturbations of the boundary can result in
    unique area minimizing surfaces that are widely separated even
    though their boundaries are nearly identical.  It was noted by
    M.~Beeson \cite{Beeson77} that a consequence of such discontinuous
    behavior is that, in a certain formal system, the area minimizing
    surface is not computable.  Thus we believe that it is essential
    to seek an approximation to an area minimizing surface the
    boundary of which is near to, but not necessarily identical to,
    the given boundary.
  \end{itemize}   
\end{ssremarks}

The last two items above concerning uniqueness and non-uniqueness
present the crucial difficulties in closely approximating the location
of an area minimizing surface, because a surface of nearly minimum
area for the given boundary may be far away in location from any area
minimizer for that boundary.  We deal with these difficulties by using
a sequence of more and more precise approximations in which we first
construct a surface \(T_{i}\)
of nearly minimum area, and then second consider an auxiliary
minimization problem.  This auxiliary problem seeks the minimum area
among surfaces satisfying two constraints which we describe informally
as follows.  The first constraint is that the boundary of each of the
surfaces considered must equal the boundary of an appropriate portion,
\(T'_{i}\),
of the surface \(T_{i}\).
The second constraint is that each of the surfaces must be relatively
far from \(T'_{i}\) in the flat norm.

Continuing our informal discussion, if \(\epsilon>0\) is specified at
the outset and if the parameters defining large and small and near and
far are chosen correctly vis-\`{a}-vis that \(\epsilon\), then in the
above sequence of constructions and minimizations, it eventually must
happen both that \(T'_{i}\) differs little from \(T_{i}\) and that the
minimum area among the surfaces considered in the auxiliary problem is
relatively large.  Consequently, the \(T'_{i}\) constructed at that
iteration is such that all surfaces relatively far from \(T'_{i}\)
have relatively large area.  Thus the surfaces with relatively small
area all must be relatively near to \(T'_{i}\).

In the previous papers \cite{Parks77} and \cite{Parks86}, the
theoretical basis was developed for computing approximations to area
minimizing surfaces by numerically approximating functions of least
gradient.  Those papers required that the given boundary for which an
area minimizing spanning surface was sought must lie on the surface of
a convex set.  An important feature of the results in those papers was
that one could be certain, at least in principle, of when sufficient
computation had been done to guarantee any desired accuracy of the
approximation in the sense of Hausdorff distance.

The method described in \cite{Parks77} and \cite{Parks86} was
implemented numerically in \cite{Parks92}.  The results reported there
and later results in \cite{Parks93} showed that, in practice, the
method gives much better approximations than the theorems
of \cite{Parks77} and \cite{Parks86} guarantee.

The requirement of \cite{Parks77} and \cite{Parks86} that the boundary
lie on the surface of a convex set is often not met.  Various
alternative methods are available for application in these
circumstances.  These are developed in the extremely general covering
space approach of K.~Brakke \cite{Brakke95}, in the duality approach
in the thesis of J.~Sullivan \cite{Sullivan90}, the more general work
of K.~Brakke \cite{Brakke95b}, and in the modification of the least
gradient method in our previous work \cite{ParksPitts96} and
\cite{ParksPitts97}.  The results of \cite{Sullivan90} and
\cite{Brakke95b} provide a way to approximate the area of the area
minimizing surface (but not the position), and implicitly so do the
results of \cite{ParksPitts97}.

We dedicate this paper to the memory of our thesis advisor and friend
Frederick~J. Almgren, Jr.

\section{The Algorithm}
\label{sec:surfaces}




The Approximation Theorem obtained by Federer and Fleming tells us
that any integral current can be approximated arbitrarily well by an
integral polyhedral chain. Consequently, given a smooth, compact,
embedded, \((N-2)\)-dimensional
boundary in \(\bbR^{N}\),
an area minimizing surface spanning the given boundary can be obtained
as the limit of integral polyhedral chains obtained by minimizing mass
in an increasing family of finite dimensional subspaces of the vector
space of \((N-1)\)-dimensional
polyhedral chains, \(\cP_{N-1}(\bbR^{N})\).
As a computational method, the obvious shortcoming of such an approach
is that, if one has in mind a desired level of accuracy of
approximation, there is no way to know whether one has achieved
it. What is lacking is \textit{a~priori} information on which finite dimension
subspace of \(\cP_{N-1}(\bbR^{N})\)
is required to obtain the desired accuracy of approximation.

In his thesis [Sul90], John Sullivan has addressed this lack of 
\textit{a~priori} information. Sullivan's approximation is carried
out using an appropriate cell complex obtained by slicing space 
with equally spaced parallel planes in each of many directions,
a structure that he calls a ``multigrid.'' 

\begin{ssdefinition}
  A {\em multigrid in} \(\bbR^{N}\)
  is the set of chains generated by a finite family of convex
  polyhedra in \(\bbR^{N}\)
  and by their vertices, edges, and faces.  In our implementation, we
  need include only the \((N-1)\)-dimensional
  faces and \((N-2)\)-faces.
\end{ssdefinition}  

Sullivan's approximation result is the following:

\begin{sstheorem}[Sul90, Theorem~6.1]\label{sullivan_main}
Given \(\epsilon\) and an
\((N-1)\)-current \(T\), we can pick a multigrid \(C\) such that \(T\) 
has a good approximation \(S\),
which is a chain in \(C\), is flat close to \(T\), and has not 
much more mass,
\(\bbM(S) \leq (1+\epsilon)\,\bbM(T)\). In fact the choice of \(C\) 
can be made merely knowing
\(\epsilon\) and bounds on \(\bbM(T)\) and on the mass of its boundary.
\end{sstheorem}

Using this last approximation result, Sullivan obtains the next result
(which we paraphrase)
regarding an algorithm for approximating the minimum area that
is required to span a given boundary cycle.

\begin{sstheorem}[Sul90, Corollary~6.2]\label{sullivan_boundary}
{\sl  Given any boundary cycle in \(\bbR^{N}\), with some a priori lower
  bound on the area of a possible area-minimizing surface, 
  a surface with no more than \(1+\epsilon\) times the true minimum area
  can be found by solving a linear programming problem.}
\end{sstheorem}

In the statement of Theorem~\ref{sullivan_boundary}, Sullivan focuses
on the approximation of the minimum area.  But we note that in
Theorem~\ref{sullivan_main} the approximating surface also
approximates the given boundary; a fact that is important in our work.
By making use of the top-dimensional polyhedra in a sequence of finer
and finer multltigrids, we are able to obtain an algorithm that not
only approximates the minimum area, but that also approximates both
the area and the location (in the sense of the \(\cF\)-norm)
of an area minimizer with boundary nearly equal to the given
boundary.  This algorithm is the first to accomplish that goal.

\begin{sstheorem}
Let \(B \in \bbI_{N-2}\) with \(\partial B = 0\) and smooth support be given. 
Let \(\epsilon>0\) be given. 
Let an open set, \(U\), with \(\spt(B)\subset U\) be given.
Then there is a computation
requiring finitely many multigrid minimizations 
that results in  a  \(T\)  guaranteed to satisfy  
the following requirements:
\begin{enumerate}
\item\label{it:contr.bdry} 
\(\spt(\partial T) \subset U\),
\item\label{it:contr.hd.bdry} 
\(\dist_{H}[ \, \spt (\partial T),\, \spt (B) \, ] < \epsilon\),
\item\label{it:contr.flt.bdry}  
\(B = \partial S + \partial T \) with
\(\spt(S) \subset U\) and \(\bbM(S) < \epsilon\),
\item\label{it:contr.mass} 
\(\bbM(T) < \epsilon +
\inf\{\bbM(S) : \partial S = B \}\),
\item\label{it:contr.flat} every area minimizing 
current \(R\) with \(\partial R = \partial T \)
is within \(\cF\)-distance \(\epsilon\) of~\(T\).
\end{enumerate}
\end{sstheorem}

\medskip\noindent{\bf Proof.}  Let \(B \in \bbI_{N-2}\)
with \(\partial B = 0\)
and smooth support be given. Let \(\epsilon>0\)
be given.  Let the open set \(U\) with \(\spt(B)\subset U\) be given.

For each \(0<r\), set 
\[
I(r) = \{ x: \dist(x,\spt B) < r \}\,,
\qquad
O(r) = \{ x: \dist(x,\spt B) \geq r \}
\,.
\]

Let \(0< \epsilon_i\), \(i=1,2,\dots\), be a decreasing sequence with limit \(0\).
Choose \(\epsilon_1\) so that
\begin{itemize}
\item
\(\epsilon_1 < \epsilon/4\),
\item
\(\hbox{\rm Clos}[ I(2\,\epsilon_1) ] \subset U\),
\item
\(\|R\|[ I(\epsilon_{1}) ] < \epsilon/3\) holds
for any mass minimizer with \(\partial R = B\),
which we can do by Proposition~5.6 of [Sul90].
\end{itemize}

\bigskip 
For each \(i\), use Sullivan's approximation method
(Theorem~\ref{sullivan_main}) to form a multigrid \(\cG(i)\) such that
for any mass minimizer \(R\) with \(\partial R = B\)
there exists \(\widehat{R}\in \cG(i)\cap \cP_2\)   such that
\begin{itemize}
\item
there exists \(S\) with
\(B = \partial S + \partial \widehat{R}\),
\(\dist_{H}[ \, \spt(S),\,\spt(B)\, ] < \epsilon_i \), 
and  \(\bbM(S) < \epsilon_i  \),
\item
\(\dist_{H}[ \, \spt (\partial \widehat{R}),\, \spt (B) \, ] < \epsilon_i \),
\item
\(\spt (\partial \widehat{R}) \subset U\),
\item
\( \bbM(\widehat{R}) \leq \bbM(R) +\epsilon_i \).
\end{itemize}
Choose the multigrids so that 
\(\cG(1)\subset \cG(2)\subset\cG(3)\subset \cdots\).

\medskip 
For each \(i\), let \(\cT(i)\subset \cG(i)\cap\cP_{N-1}\) be the set of 
currents, \(T\), satisfying
\begin{itemize}
\item
there exists \(S\) with
\(B = \partial S + \partial T\),
\(\dist_{H}[ \, \spt(S),\,\spt(B)\, ] < \epsilon_i\), and 
\(\bbM(S) < \epsilon_i \),

\item
\(\dist_{H}[ \, \spt (\partial T),\, \spt (B) \, ] < \epsilon_i\),
\item
\(\spt (\partial T) \subset U\).
\end{itemize}

\noindent
Using an appropriate algorithm, obtain \(T_i\in \cT(i)\) such that
$$
\bbM(T_i) \leq \epsilon_i + \inf\{\bbM(T):T\in \cT(i) \}
\,.
$$
(We are solving a linear programming problem. We are also not requiring the
exact solution; only that we be within \(\epsilon_{i}\) of the minimum
value of the objective function.)

\bigskip\noindent{\bf Claim 1.} If \(\mu\) denotes
the mass of any mass minimizer \(R\) with \(\partial R = B\),
then
\begin{equation}\label{(*)}
\mu - \epsilon_i \leq \bbM(T_i) \leq \mu + 2\,\epsilon_i 
\end{equation}
holds  for each \(i\), 
and the limit of any  \(\cF\)-convergent 
subsequence of 
\(\big\{ T_i \big\}_{i = 1}^\infty\)
is a mass minimizer with boundary equal to \(B\).

\medskip\noindent
{\bf Proof of Claim.}
Let \(R\) be a mass minimizer with \(\partial R = B\). 

Since \(T_i\in \cT(i)\), there exists \(S_i\) with 
\(B = \partial S_i + \partial T_i = \partial (S_i+T_i) \) 
and \(\bbM(S_i) < \epsilon_i \).
Thus we have
$$
\mu = \bbM(R) \leq \bbM(S_i+T_i) \leq \bbM(S_i) + \bbM(T_i)
\leq \epsilon_i +\bbM(T_i)
\,,
$$
giving us the left-hand inequality in (\ref{(*)}).

We have chosen the multigrid \({\cG}(i)\)
so that for any mass minimizer \(R\)
with \(\partial R = B\)
there exists \(\widehat{R}\in \cT(i)\) such that
\begin{itemize}
\item
there exists \(S\) with
\(B = \partial S + \partial \widehat{R}\),
\(\dist_{H}[ \, \spt(S),\,\spt(B)\, ] < \epsilon_i\), 
and  \(\bbM(S) < \epsilon_i \),
\item
\(\dist_{H}[ \, \spt (\partial \widehat{R}),\, \spt (B) \, ] < \epsilon_i\),
\item
\(\spt (\partial \widehat{R}) \subset U\),
\item
\( \bbM(\widehat{R}) \leq \bbM(R) +\epsilon_i\).
\end{itemize}
Then \(\widehat{R}\) satisfies the conditions for membership in \(\cT(i)\).
By the choice of \(T_i\), we conclude that 
$$
\bbM(T_i) \leq \epsilon_i + \bbM(\widehat{R}) 
\leq 2\,\epsilon_i + \bbM(R) = 2\,\epsilon_i + \mu
\,,
$$
giving us the right-hand inequality in  (\ref{(*)}).

\medskip
Now, let \(T^*\) be the limit of any  \(\cF\)-convergent 
subsequence of 
\(\big\{ T_i \big\}_{i = 1}^\infty\).
Passing to that subsequence, but without changing notation,
we suppose \(T_i\ra T^*\). 
Letting \(S_i\) be as above, we have \(B= \partial(S_i+T_i)\)
and \(\bbM(S_i)\ra 0\). So \(B = \partial T^*\).
By the lower semicontinuity of  mass, 
\(\bbM(T^*)\leq \lim_{i\ra\infty} \bbM(T_i) = \mu\).
Thus
\(T^*\)
is a mass minimizer with boundary \(B\). 

\noindent{\bf Claim 1 has been proved.}

\bigskip\noindent{\bf Claim 2.} For infinitely many \(i\), we 
have
$$
\bbM[ T_i\elbow I(\epsilon_{1}) ] \leq \epsilon/2
\,.
$$

\medskip\noindent
{\bf Proof of Claim.} Suppose Claim~2 were false. Then
there would be but finitely many elements in
$$
J = \{ i : \bbM[ T_i\elbow I(\epsilon_{1}) ] \leq \epsilon/2 \}
\,.
$$
Set  \(i_0= 1+ \max J\). Then 
$$
\bbM[ T_i\elbow I(\epsilon_{1}) ] > \epsilon/2
$$
holds for all \(i \geq i_0\).
Since  
$$
\bbM[T_i] = \bbM[ T_i\elbow I(\epsilon_{1}) ]
+ \bbM[ T_i\elbow O(\epsilon_1) ]
$$
we have
$$
\bbM[ T_i\elbow O(\epsilon_1) ]
= \bbM[T_i] - \bbM[ T_i\elbow I(\epsilon_{1}) ]
< \bbM[T_i] - \epsilon/2
\,.
$$
So
$$
\lim_{i\ra\infty} \bbM[ T_i\elbow O(\epsilon_1) ]
\leq \mu -\epsilon/2,
$$
where, as in Claim~1, \(\mu\) denotes the mass of any minimizer with
boundary~\(B\).

Passing to an \(\cF\)-convergent
subsequence, but without changing notation, we 
may suppose \(T_i\) converges to 
a mass minimizer \(R\) with \(\partial R = B\).
By the lower semicontinuity of  mass,
$$
\| R\| [O(\epsilon_1)] \leq  \mu -\epsilon/2
$$
holds.
Since \(\bbM[R] = \mu\), we have
$$
\| R\| [I(\epsilon_{1})]  \geq \epsilon/2
\,,
$$ 
contradicting the requirement in the definition of \(\epsilon_1\) that
\(\|R\|[ I(\epsilon_{1}) ] < \epsilon/3\) hold.

\noindent{\bf Claim 2 has been proved.}

\bigskip
Let \(\cK\) be a closed set disjoint from
\(I(\epsilon_1/2)\),
containing \(O(\epsilon_1)\),  and having a  
polyhedral boundary.
For each \(i=1,2,\dots\), set 
$$
T'_i =  T_{i}\elbow \cK
\hbox{\rm\quad and\quad}
B_i = \partial T'_i
\,.
$$ 

For each \(i\),
use Sullivan's approximation method (Theorem~\ref{sullivan_main}) to form a
multigrid \(\cG'(i)\),
with \(\cG(i) \subset \cG'(i)\)
and \(T'_i\in \cG'(i)\),
such that for any mass minimizer \(R\)
with \(\partial R = B_i\)
there exists \(\widehat{R}\in \cG(i)\cap \cP_{N-1}\) such that
\begin{itemize}
\item
there exists \(S\in \cG'(i) \cap \cP_{N-1}\) with
\(B_{i } = \partial S  + \partial \widehat{R}\),
\(\bbM(S ) < \epsilon_{i} \), and
\(\dist_{H}[ \, \spt(S ),\,\spt(B_{i })\, ] < \epsilon_{i }\),
\item
\(\dist_{H}[ \, \spt (\partial \widehat{R}),
\, \spt (B_{i}) \, ] < \epsilon_{i}\),
\item
\(\spt (\partial \widehat{R}) \subset U\),
\item
\( \bbM(\widehat{R}) \leq \bbM(R) +\epsilon_{i}\),
\item
\(\widehat{R} - R = X + \partial Y\) for some \(X\) and \(Y\) with
\(
\bbM(X) + \bbM(Y)  \leq \epsilon_{i }
\).
\end{itemize}
Choose the multigrids so that
\(\cG'(1)\subset \cG'(2)\subset\cG'(3)\subset \cdots\).

\medskip
For each \(i\), 
let  \(\cT'(i)\subset \cG'(i)\cap \cP_{N-1}\) be the set of  currents,
\(T\), satisfying
\begin{itemize}
\item
there exists \(S\in \cG'(i) \cap \cP_{N-1}\) with
\(B_i = \partial S + \partial T\),
\(\bbM(S) < \epsilon_{i} \),
and  \(\dist_{H}[ \, \spt(S),\,\spt(B_i)\, ] < \epsilon_{i}\),
\item
\(\dist_{H}[ \, \spt (\partial T),\, \spt (B_i) \, ] < \epsilon_{i}\),
\item
\(\spt (\partial T) \subset U\).
\end{itemize}

\noindent
For each \(i\), let \(\cQ(i)\subset \cT'(i)\) be the set of currents,
\(Q\), satisfying
\begin{itemize}
\item
\( \bbM(W)  \geq \epsilon/2\), where \(\partial W = T'_i-S-Q\) where 
\(S\) is as in the first condition for membership of \(Q\) in \(\cT'(i)\).
\end{itemize}
Notice that if \(B_i = \partial S + \partial Q\), then
\(W\) satisfying  \(\partial W = T'_i-S-Q\) is unique and
\(W\in \cG'(i) \cap \cP_{N}\).

\medskip\noindent
Using an appropriate algorithm, obtain \(Q_i\in \cQ(i)\) such that
$$
\bbM(Q_i) \leq \epsilon_{i} + \inf\{\bbM(Q):Q\in \cQ(i) \}.
$$
(We are solving a linear programming problem. We are also not
requiring the exact solution, only that we be within \(\epsilon_{i}\)
of the minimum value of the objective function.)

\bigskip\bigskip\noindent
{\bf Stopping Conditions:}
\begin{enumerate}[(C1)]
\item\label{eq:first.stop}
  \(\bbM(Q_{i}) \geq \bbM(T'_{i }) + 3\,\epsilon_{i}\)
\item\label{eq:second.stop}
  \(\bbM(T\elbow\bbR^{N}\setminus\cK)\leq\epsilon/2\)
\end{enumerate}


\bigskip\noindent{\bf Claim 3.}
If for some \(i_0\),  the stopping conditions
are satisfied, then \( T'_{i_0} \) is the 
desired approximation. That is,
\begin{itemize}
\item  
\(B = \partial S + \partial  T'_{i_0} \) with
\(\spt(S) \subset U\) and \(\bbM(S) < \epsilon\),
\item 
\(\dist_{H}[ \, \spt (\partial  T'_{i_0} ),\, \spt (B) \, ] < \epsilon\),
\item 
\(\spt(\partial  T'_{i_0} ) \subset U\),
\item 
\(\bbM(  T'_{i_0} ) < \epsilon +
\inf\{\bbM(S) : \partial S = B \}\),
\item  
every 
mass minimizing 
current \(R\) with \(\partial R = \partial T'_{i_0} = B_{i_0} \)
is within \(\cF\)-distance \(\epsilon\) of
\(  T'_{i_0} \).
\end{itemize}

\medskip\noindent
{\bf Proof of Claim.}
By the choice of \(\epsilon_1\), it is immediate  that
\begin{itemize}
\item 
\(\dist_{H}[ \, \spt (\partial  T'_{i_0} ),\, \spt (B) \, ] < \epsilon\) ,
\item 
\(\spt(\partial   T'_{i_0} ) \subset U\)
\end{itemize}
hold.

Since \(T_{i_0} \in \cT(i_0)\),
there exists \(S_1\) with
$$
B = \partial S_1 + \partial T_{i_0},
\ \ \dist_{H}[ \, \spt(S_1),\,\spt(B)\, ] <
\epsilon_{i_0}, \hbox{\rm\ \ and\ \ } \bbM(S_1) < \epsilon_{i_0} 
\,.
$$
So
\begin{eqnarray*}
B &=& \partial S_1 + \partial T_{i_0} \\
&=& \partial S_1 + \partial \Big(T_{i_0} \elbow  \bbR^{N}\setminus\cK
+ T_{i_0} \elbow  \cK \Big)\\
&=& \partial  \Big( S_1 + 
T_{i_0} \elbow  \bbR^{N}\setminus\cK  \Big)
+ \partial T'_{i_0}
\,.
\end{eqnarray*}
We have
$$
\spt( S_1 + 
T_{i_0} \elbow  \bbR^{N}\setminus\cK ) \subset U
$$
and
$$
\bbM\Big( S_1 + 
T_{i_0} \elbow  \bbR^{N}\setminus\cK  \Big)
\leq 
\bbM( S_1) + 
\bbM( T_{i_0} \elbow  \bbR^{N}\setminus\cK  )
\leq \epsilon_{i_0} +   \epsilon/2
\leq \epsilon
\,,
$$
where we have used the stopping condition
(C\ref{eq:second.stop}).

The right-hand inequality in (\ref{(*)}) gives us
$$
\bbM(  T'_{i_0} ) < \epsilon +
\inf\{\bbM(S) : \partial S = B \}
\,.
$$

\medskip
Suppose \(R\) is a minimizer with \(\partial R = B_{i_0}\).
Let  \(\widehat{R}\) be such that
\begin{itemize}
\item
there exists \(S_2\) with
\(B_{i_0} = \partial S_2 + \partial \widehat{R}\),
\(\dist_{H}[ \, \spt(S_2),\,\spt(B_{i_0})\, ] < \epsilon_{i_0}\), 
and  \(\bbM(S_2) < \epsilon_{i_0} \),
\item
\(\dist_{H}[ \, \spt (\partial \widehat{R}),
\, \spt (B_{i_0}) \, ] < \epsilon_{i_0}\),
\item
\(\spt (\partial \widehat{R}) \subset U\),
\item
\( \bbM(\widehat{R}) \leq \bbM(R) +\epsilon_{i_0}\),
\item
\(\widehat{R} - R = X + \partial Y\) for some \(X\) and \(Y\) with
\(
\bbM(X) + \bbM(Y)  \leq \epsilon_{i_0}
\).
\end{itemize}
Notice that the first three conditions above tell us that
\(\widehat{R}\in \cT'(i_0)\).

Next, note that since \(R\) is a mass minimizer with 
\(\partial R = \partial T'_{i_0}\), we have
$$
\bbM(R) \leq  \bbM( T'_{i_0}) 
\,.
$$
Thus 
$$
\bbM(\widehat{R}) \leq  \bbM(R) + \epsilon_{i_0}
\leq 
\bbM(T'_{i_0}) + \epsilon_{i_0}
$$
holds.
If it were the case that \(\widehat{R}\in \cQ(i_0)\),
then the choice of \(Q_{i_0}\) would give us
$$
\bbM(Q_{i_0}) \leq \epsilon_{i_0} + \bbM(\widehat{R})
\leq 
\bbM(T'_{i_0}) + 2\,\epsilon_{i_0}
\,,
$$
contradicting the stopping condition (C\ref{eq:first.stop}).
We conclude that \(\widehat{R}\in \cT'(i_0)\setminus\cQ(i_0)\).

Now, let  \(W\) satisfy \(\partial W = T'_{i_{0}} - S_{2} - \widehat{R}\) with 
\(S_2\) as above.
Because \(\widehat{R} \notin \cQ(i_0)\), we have
$$ 
\bbM(W)  < \epsilon/2
\,.
$$

We also have 
\(\widehat{R} - R = X + \partial Y\) for some \(X\) and \(Y\) with
$$
\bbM(X) + \bbM(Y)  \leq \epsilon_{i_0}
\,.
$$
Consequently, we see that 
$$
T'_{i_0} - R =
S_2 + X + \partial Y + \partial W
\,,
$$
with
$$
\bbM(S_2) + \bbM(X) + \bbM(Y) +\bbM(W)
\leq 2\, \epsilon_{i_0} + \epsilon/2 \leq \epsilon
\,.
$$
That is, we have \(\cF(T'_{i_0} - R)  \leq \epsilon\).

\noindent{\bf Claim 3 has been proved.}

\bigskip\noindent{\bf Claim 4.} For some \(i\), the stopping conditions
will be satisfied.

\medskip\noindent {\bf Proof of Claim.} Applying Claim~2, we pass to a
subsequence (without changing notation) for which the stopping
condition~(C\ref{eq:second.stop}) holds for all \(i\).

Arguing by contradiction, we suppose that
$$
\bbM(Q_{i}) <  \bbM(T'_{i}) + 3\,\epsilon_{i}
$$
holds for every \(i\).

Since  \(T_{i }\in \cT(i )\), there exists   \(S_{i }\) with
\(B = \partial S_{i } + \partial T_{i } \)
and \(\bbM(S_{i } ) < \epsilon_{i }\).
Since \(Q_i\in \cQ(i)\), there exists 
\(S'_i\) with 
\(\partial  T'_{i} = B_i = \partial S'_i + \partial Q_i\)
and \(\bbM(S'_i) < \epsilon_{i}\).

Set
$$
P_i = S_{i } +  T_{i }\elbow \bbR^{N}\setminus\cK + S'_i+ Q_i
\,.
$$
We have 
\begin{eqnarray*}
\partial P_i &=& \partial S_i +
\partial [T_{i }\elbow \bbR^{N}\setminus\cK]
+ \partial S'_i + \partial Q_i\\
&=&\partial S_i +
\partial [T_{i }\elbow \bbR^{N}\setminus\cK] + \partial T'_i\\
&=& \partial S_i + \partial T_i 
\ = \ B
\end{eqnarray*}
and
\begin{eqnarray*}
\bbM(P_i) &\leq &
\bbM( S_{i } ) +  \bbM[ T_{i }\elbow \bbR^{N}\setminus\cK ]+ 
\bbM(S'_i) + \bbM( Q_i )\\
&\leq & 2\,\epsilon_i + \bbM[ T_{i }\elbow \bbR^{N}\setminus\cK ]+
\bbM(T'_i) + 3\,\epsilon_i
\ =\ \bbM(T_i) + 5\,\epsilon_i
\,.
\end{eqnarray*}

We may pass to a subsequence,
again without changing notation, such
that \(P_i\) converges to \(P^*\) and \(S_i + T_{i }\) converges \(T^*\).
By the lower semicontinuity of mass and the right-hand 
inequality in (\ref{(*)}), we see that both
\(P^*\) and \( T^*\) are  mass
minimizers with boundary \( B\). By construction, \(P^*\) and \( T^*\)
are equal in \(I(\epsilon_1/2)\).  By the regularity theory of mass
minimizers, the singular set of a minimizer cannot disconnect the 
surface.  We have \(P^*= T^*\).

The fact that \(P^* = T^*\) tells us that 
\(\cF[ (S_i+ T_i) - P_{i} ] \ra 0\), so 
we can write   \((S_i+ T_i) - P_{i} = X_i + \partial Y_i\)
with \(\bbM(X_i) + \bbM(Y_i) \ra 0\). Then applying the isoperimetric 
inequality to \(X_i\), we see that we can write 
\((S_i+ T_i) - P_{i}  =   \partial Z_i\) with \(\bbM(Z_i)\ra 0\).

On the other hand, observe that
$$
(S_i+ T_i) - P_{i} = T'_i -Q_i - S'_i
\,.
$$
By the definition of \(\cQ(i)\), we have
\(T'_i -Q_i - S'_i = \partial W_i\) with \(\bbM(W_i) \geq \epsilon/2\).
This last inequality contradicts   \(\bbM(Z_i)\ra 0\),
because \(W_i\) and \(Z_i\) are \(N\)-dimensional integral currents
in \(\bbR^N\) having the same boundary, so in fact, they are equal.

\noindent{\bf Claim 4 has been proved.}

\bigskip
\noindent{\bf Conclusion.}
Once the sequence \(\epsilon_i\) satisfying the required
conditions has been chosen, 
the algorithm proceeds as follows:
\begin{enumerate}
\item[(A1)] Set \(i = 1\).
\item[(A2)]
Compute \(T_i\).
\item[(A3)]
If  the condition 
\(\bbM[ T_{i }\elbow \bbR^{N}\setminus\cK ] \leq \epsilon/2 \)
is satisfied, then advance to step (A4). 
Otherwise, increment \(i\) and go to step (A2).
\item[(A4)]
Compute \(Q_{i}\).
\item[(A5)]
If the  condition 
\(\bbM(Q_{i}) \geq  \bbM( T_{i}\elbow \cK) + 3\,\epsilon_{i}\)
is satisfied, then return \(T'_i \) and terminate the algorithm. 
Otherwise, increment  \(i\) and go to step (A2).
\end{enumerate}

Claim~4 guarantees that the algorithm terminates after finitely 
many steps, while Claim~3 guarantees that the returned 
value \(T'_{i} \) is the desired approximation.
\bibliography{main}
\bibliographystyle{amsalpha}

\end{document}